\input amstex
\documentstyle{amsppt}
\topmatter
\title Poincar\'{e} series and modular functions for $U(n, 1)$\endtitle
\author {\rm LEI YANG} \endauthor
\abstract \nofrills \it \qquad \qquad \qquad \qquad \quad Peking University,
Beijing, P. R. China \endabstract

\address
\flushpar \qquad {\smc Lei Yang}
\flushpar \qquad {\rm Department of Mathematics}
\flushpar \qquad {\rm Peking University}
\flushpar \qquad {\rm Beijing, 100871}
\flushpar \qquad {\rm P. R. China}
\flushpar \qquad {\rm yanglei\@sxx0.math.pku.edu.cn}
\endaddress
\endtopmatter
\document

\centerline{{\bf 1. Introduction}}
\vskip 0.5 cm
\flushpar
  In the theory of automorphic forms, two classes of rank one
reductive Lie groups $O(n, 1)$ and $U(n, 1)$ are the important
objects. Automorphic forms on $O(n, 1)$ have been intensively
studied. In this paper we study the automorphic forms on
$U(n, 1)$. We construct infinitely many modular forms and
non-holomorphic automorphic forms on $U(n, 1)$ with respect to
a discrete subgroup of infinite covolume. More precisely,
we obtain the following theorem:

\proclaim\nofrills{\smc Main Theorem}. {\it A function
$f: {\Bbb H}_{\bold C}^{n+1}:=\{ Z=(z_1, \cdots, z_{n+1}) \in
{\bold C}^{n+1}: \text{Im}z_{n+1}>\sum_{j=1}^{n} |z_j|^2 \}
\to {\bold C}$ is called a nonholomorphic automorphic form attached to the
unitary group $U(n+1, 1)$ if it satisfies the following three
conditions:
\roster
\item
$f$ is an eigenfunction of the Laplace-Beltrami operator
$L$ of $U(n+1, 1)$ on ${\Bbb H}_{\bold C}^{n+1}$;
\item
$f$ is invariant under the modualr group; i.e.,
$f(\gamma(Z))=f(Z) $ for all $\gamma \in G({\bold Z})$ and
all $Z \in {\Bbb H}_{\bold C}^{n+1}$, where $G({\bold Z})$ is
the discrete subgroup of $G:=TU(n+1, 1)T^{-1}$
defined over ${\bold Z}$ and $T$ is the Cayley transform
from the unit ball ${\bold B}_{\bold C}^{n+1}$ in
${\bold C}^{n+1}$ to ${\Bbb H}_{\bold C}^{n+1}$.
\item
$f$ has at most polynomial growth at infinity; i.e., there
are constants $C>0$ and $k$ such that $|f(Z)| \leq C \rho^k$,
as $\rho \to \infty$ uniformly in $t$, for fixed $\beta$.
Here $\rho(Z)=\text{Im}z_{n+1}-\sum_{j=1}^{n} |z_j|^2$ and
$\beta(Z)=\sum_{j=1}^{n} |z_j|^2$.
\endroster

We denote by ${\Cal N}(G({\bold Z}), \lambda)$
the space of such nonholomorphic automorphic forms
attached to $U(n+1, 1)$. Then there exist a family
of Poincar\'{e} series (hence, infinitely many
elements)
$r(Z, Z^{\prime}; s) \in {\Cal N}(G({\bold Z}),
s(s-n-1))$ for $\text{Re}(s)>n$, where
$$r(Z, Z^{\prime}; s):=\sum_{\gamma \in
  G({\bold Z})} g_s(x(Z, \gamma(Z^{\prime})),
  y(Z)). \eqno(1.1)$$
In the nondegenerate case,
$$g_s(x, y)=g_s(x, y; a, b)=x^{-a} y^{-b}
  F_3(a, b; a, b-n+1; 2s-n; -x^{-1}, -y^{-1})
  \eqno(1.2)$$
for $\text{Re}(a)>1$, $\text{Re}(b)>n-1$ with
$a+b=s$ and $F_3=F_3(\alpha, \alpha^{\prime}; \beta,
\beta^{\prime}; \gamma; x, y)$ is a two variable
hypergeometric function; in the degenerate case,
$$g_s(x, y)=\left\{\aligned
  &x^{-s} {}_2F_1(s, s; 2s-n; -x^{-1}), \quad
   \text{or}\\
  &(x+y)^{-s} {}_2F_1(s, s-n; 2s-n; -(x+y)^{-1}),
\endaligned\right.\eqno(1.3)$$
where ${}_2F_1={}_2F_1(\alpha, \beta; \gamma; z)$
is a hypergeometric function. Here
$$x(Z, Z^{\prime})=\frac{(t(Z)-t(Z^{\prime}))^2+
  (\rho(Z)+\beta(Z)-\rho(Z^{\prime})-\beta(
  Z^{\prime}))^2}{4 \rho(Z) (\rho(Z^{\prime})+
  \beta(Z^{\prime}))}, \quad
  y(Z)=\frac{\beta(Z)}{\rho(Z)} \eqno(1.4)$$
are two invariants under the action of
$G({\bold Z})$, where $Z, Z^{\prime} \in
{\Bbb H}_{\bold C}^{n+1}$ and $t(Z)=\text{Re}
(z_{n+1})$.

  On the other hand, the Eisenstein series
$E(Z; \lambda, \mu) \in {\Cal N}(G({\bold Z}), s(s-n-1))$,
where
$$E(Z; \lambda, \mu)=\sum_{\gamma \in G({\bold Z})_{\infty}
  \backslash G({\bold Z})} \rho(\gamma(Z))^{\lambda}
  \beta(\gamma(Z))^{\mu} \eqno (1.5)$$
with $(\lambda, \mu)=(s, 0)$ or $(s, 1-n)$.

  A modular form on $U(n+1, 1)$ associated
with $G({\bold Z})$ is a function $\phi:
{\Bbb H}_{\bold C}^{n+1} \to {\bold C}$ satisfying
the following transform equations:
\roster
\item
$\phi(\frac{z}{cz_{n+1}+d}, \frac{az_{n+1}+b}{cz_{n+1}+d})
=(c z_{n+1}+d)^k e^{2 \pi i m c(z_1^2+\cdots+z_n^2)/(c z_{n+1}+d)}
\phi(z, z_{n+1})$, for $z=(z_1, \cdots, z_n)$.
\item
$\phi(wz, z_{n+1})=\phi(z, z_{n+1})$
for all $w \in S_n$, where $S_n$ is the symmetric group of $n$-order.
\item
$\phi(z, z_{n+1})$ is a locally bounded function as
$\text{Im} z_{n+1} \to \infty$.
\endroster

  The set of modular forms on $U(n+1, 1)$ is denoted as
$M_{k, m}(G({\bold Z}))$. Then
$$j_m(z, z_{n+1}) :=\frac{1728 g_{2, m}(z, z_{n+1})^3}
 {\Delta_m(z, z_{n+1})} \in M_{0, 0}(G({\bold Z})),
 \eqno(1.6)$$
where
$g_{2, m_1}(z, z_{n+1}) :=\frac{4}{3} \pi^4
  E_{4, m_1}(z, z_{n+1})$,
$g_{3, m_2}(z, z_{n+1}) :=\frac{8}{27} \pi^6
  E_{6, m_2}(z, z_{n+1})$,
and
$\Delta_m(z, z_{n+1}) :=g_{2, m}(z, z_{n+1})^3-27
  g_{3, \frac{3}{2}m}(z, z_{n+1})^2.$
Moreover, $\{ j_m \}_{m \in {\bold N}}$ are a family of
modular functions on the modular variety ${\Cal M}_{n+1}
:=G({\bold Z}) \backslash {\Bbb H}_{\bold C}^{n+1}$.}
\endproclaim

\vskip 0.5 cm
\centerline{\bf 2. The Laplace-Beltrami operator and the discrete}
\centerline{\bf subgroup of $U(n+1, 1)$ on ${\Bbb H}_{\bold C}^{n+1}$}
\vskip 0.5 cm

\flushpar
  For $g=\left(\matrix
                A & B\\
                C & D
               \endmatrix\right) \in U(n+1, 1)$,
where $A=(a_{ij})_{(n+1) \times (n+1)}$,
$B=(b_{ij})_{(n+1) \times 1}$,
$C=(c_{ij})_{1 \times (n+1)}$, $D=(d)_{1 \times 1}$, and
$W=(w_1, \cdots, w_{n+1}) \in {\bold B}_{\bold C}^{n+1}
=\{ W \in {\bold C}^{n+1} | |W|^2=|w_1|^2+\cdots+|w_{n+1}|^2
< 1 \}$, the action of $g$ on $W$ is defined as
$g \circ W=[(A W^{t}+B) (C W^{t}+D)^{-1}]^{t}$.
According to \cite{5}, the corresponding
Laplace-Beltrami operator on $U(n+1, 1)$ is
$$L_{{\bold B}^{n+1}}=\text{tr}[(I-W W^{*})
  \overline{\partial_W} \cdot (I-W^{*} W) \cdot
  \partial_{W}^{\prime}], \eqno (2.1)$$
where $\partial_W$ is the differential operator
$\partial_W=\left( \frac{\partial}{\partial w_1},
\cdots, \frac{\partial}{\partial w_{n+1}} \right),$
and the dots here indicate that the factor
$(I-W^{*} W)$ is not differentiated. 

\proclaim\nofrills{\smc Theorem 2.1}. {\it The Laplace-Beltrami operator
of $U(n+1, 1)$ on ${\bold B}_{\bold C}^{n+1}$ is
$$\Delta=(1-\sum_{j=1}^{n+1} |w_j|^2)
  \left[\sum_{j=1}^{n+1} \frac{\partial^2}{\partial w_j
  \partial \overline{w_j}}-\left(\sum_{j=1}^{n+1} w_j
  \frac{\partial} {\partial w_j} \right)
  \left(\sum_{j=1}^{n+1} \overline{w_j}
  \frac{\partial}{\partial \overline{w_j}} \right)
  \right]. \eqno (2.2)$$}
\endproclaim

  Let ${\Bbb H}_{\bold C}^{n+1}$ be the Siegel domain
of the second kind, i.e., the complex hyperbolic space,
${\Bbb H}_{\bold C}^{n+1}=\{ (z_1, \cdots, z_{n+1}) \in
  {\bold C}^{n+1} | \text{Im} z_{n+1} >
  \sum_{j=1}^{n} |z_j|^2 \}.$
It is well-known that ${\Bbb H}_{\bold C}^{n+1}$ is
holomorphically equivalent to ${\bold B}_{\bold C}^{n+1}$.
The Cayley transform $T$ from ${\bold B}_{\bold C}^{n+1}$
onto ${\Bbb H}_{\bold C}^{n+1}$ is given by
$z_1=\frac{i w_1}{1-w_{n+1}}, \cdots,
 z_n=\frac{i w_n}{1-w_{n+1}},
 z_{n+1}=i \frac{1+w_{n+1}}{1-w_{n+1}}.$

\proclaim\nofrills{\smc Theorem 2.2}. {\it The Laplace-Beltrami
operator of $U(n+1, 1)$ on the complex hyperbolic space
${\Bbb H}_{\bold C}^{n+1}$ is
$$L=(\text{Im}z_{n+1}-\sum_{j=1}^{n} |z_j|^2)
  \left[ \sum_{j=1}^{n+1} 2i \left(\overline{z_j}
  \frac{\partial^2}{\partial z_{n+1} \partial \overline{z_j}}
  -z_j \frac{\partial^2}{\partial \overline{z_{n+1}}
  \partial z_j} \right)+\sum_{j=1}^n \frac{\partial^2}
  {\partial z_j \partial \overline{z_j}} \right].
  \eqno (2.3)$$}
\endproclaim

  Let us consider ${\bold B}_{\bold C}^{n+1}$ as a
symmetric space $U(n+1, 1)/U(n+1) \times U(1)$. Set 
$G=T U(n+1, 1) T^{-1}$. Let $K$ be a maximal compact
subgroup of $G$ and $\Gamma$ be a discrete subgroup of
$G$. ${\Bbb H}_{\bold C}^{n+1}=G/K$ is invariant under
the action of $G$.

  Let $H$ be a subgroup of $U(n+1, 1)$, if $T H T^{-1}=\Gamma$
is a discrete subgroup of $G$, then
$$\left(\matrix
   i I_n &    &  \\
         &  i & i\\
         & -1 & 1
\endmatrix\right) \left(\matrix
     A & B\\
     C & D
\endmatrix\right)=\left(\matrix
    \widetilde{A} & \widetilde{B}\\
    \widetilde{C} & \widetilde{D} 
\endmatrix\right) \left(\matrix
    i I_n &    &   \\
          &  i & i \\
          & -1 & 1
\endmatrix\right),$$
where $\left(\matrix
              A & B\\
              C & D
        \endmatrix\right) \in H$ and
$\left(\matrix
       \widetilde{A} & \widetilde{B}\\
       \widetilde{C} & \widetilde{D}
 \endmatrix\right) \in \Gamma$.  
We have
$\left(\matrix
        iA & iB\\
        SC & SD
       \endmatrix\right)=\left(\matrix
        i \widetilde{A} & \widetilde{B} S\\
        i \widetilde{C} & \widetilde{D} S
       \endmatrix\right),$
where $S=\left(\matrix
                i & i\\
               -1 & 1
              \endmatrix\right)$.
Thus,
$A=\widetilde{A}$, $B=-i \widetilde{B} S$,
$C=i S^{-1} \widetilde{C}$, $D=S^{-1}
\widetilde{D} S$.
By $\left(\matrix
           A & B\\
           C & D
\endmatrix\right) \in U(n+1, 1)$, we have
$$\widetilde{A} \widetilde{A}^{*}+\widetilde{B}
  J \widetilde{B}^{*}=I_n, \quad
  \widetilde{C} \widetilde{A}^{*}+\widetilde{D}
  J \widetilde{B}^{*}=0, \quad
  \widetilde{C} \widetilde{C}^{*}+\widetilde{D}
  J \widetilde{D}^{*}=J,$$
where $J=\left(\matrix
          0 & -2i\\
         2i &   0
  \endmatrix\right)$.
If $\Gamma$ is defined over ${\bold Z}$, then
$\widetilde{B}=\widetilde{C}=0$, 
$\widetilde{A} \in O(n, {\bold Z})$,
$\widetilde{D} \in SL(2, {\bold Z})$.
In fact, $\widetilde{A}$ is a permutation
matrix with elements $0$ and $\pm 1$. For simplicity,
we only consider the permutation matrix with elements
$0$ and $1$. We denote $\sigma$ as the above permutation
matrix, it can be identified with the element of the
symmetric group of $n$-order $S_n$. Therefore,
$$\Gamma=\left\{ \gamma=\left(\matrix
  \sigma &   &  \\
         & a & b\\
         & c & d
 \endmatrix\right):
 \left(\matrix
  a & b\\
  c & d
 \endmatrix\right) \in SL(2, {\bold Z})
 \right\} \eqno (2.4)$$
where $\sigma \in S_n$ is a permutation on the set
$\left\{ z_1, \cdots, z_n \right\}$ as $\sigma(z_j)=
z_{\sigma(j)}$. Thus
$$\gamma(z_1, \cdots, z_{n+1})=\left(
  \frac{\sigma(z_1)}{cz_{n+1}+d}, \cdots,
  \frac{\sigma(z_n)}{cz_{n+1}+d},
  \frac{a z_{n+1}+b}{cz_{n+1}+d} \right).
  \eqno (2.5)$$
We denote this group as $G({\bold Z})$.
Its parabolic subgroup is
$$G({\bold Z})_{\infty}=
  \left\{ \left(\matrix
  \sigma &       &      \\
         & \pm 1 & n    \\
         &     0 & \pm 1
\endmatrix \right): n \in {\bold N}, \sigma
\in S_n \right\}.$$

  In coordinates $(t, \rho, z_1, \cdots, z_n)$
where $t=\text{Re}z_{n+1}$, $\rho=\text{Im}z_{n+1}-
\sum_{j=1}^{n} |z_j|^2$, we define
$\rho(Z)=\text{Im}z_{n+1}-\sum_{j=1}^{n} |z_j|^2$ and
$\beta(Z)=\sum_{j=1}^{n} |z_j|^2$, then
$\rho \circ \gamma(Z)=\frac{\rho(Z)}{|cz_{n+1}+d|^2}$
and
$\beta \circ \gamma(Z)=\frac{\beta(Z)}{|cz_{n+1}+d|^2}$.

  Denote $\Theta$ the inverse conjugate transpose, it
is an automorphism of the Lie groups called the Cartan
involution. Let
$U(n+1, 1)=K_U A_U N_U=\overline{N_U} A_U K_U$
be the Iwasawa decomposition, where
$$N_U=\left\{ n(z, t)=\left(\matrix
  I_n & iz^{t} & -i z^{t}\\
  i \overline{z} & 1-\frac{|z|^2-it}{2} &
  \frac{|z|^2-it}{2}\\
  i \overline{z} & -\frac{|z|^2-it}{2}
  &1+\frac{|z|^2-it}{2}
\endmatrix\right): z \in {\bold C}^n, t \in
{\bold R} \right\},$$
$$\overline{N_U}=\Theta N_U,$$
$$A_U=\left\{ a(\zeta)=\left(\matrix
  I_n & 0 & 0\\
    0 & \text{ch} \zeta & \text{sh} \zeta\\
    0 & \text{sh} \zeta & \text{ch} \zeta
 \endmatrix\right): \zeta \in {\bold R} \right\},$$
$$K_U=U(n+1) \times U(1).$$
Let $P_U=\overline{N_U} A_U$ be the semidirect
product of $\overline{N_U}$ and $A_U$, where the action
of $A_U$ on $\overline{N_U}$ is given by
$$a(\zeta): \overline{n}(z, t) \mapsto a(\zeta)^{-1} 
\overline{n}(z, t) a(\zeta)=\overline{n}(e^{\zeta} z,
e^{2 \zeta} t).$$
In coordinates $z=(z_1, \cdots, z_n)$, $t=\text{Re}z_{n+1}$
and $\rho=\text{Im}z_{n+1}-|z|^2$,
${\Bbb H}_{\bold C}^{n+1}=\{ (z, t, \rho):
z \in {\bold C}^n, t \in {\bold R}, \rho >0 \}$,
$\partial {\Bbb H}_{\bold C}^{n+1}=\{ (z, t)=(z, t, 0):
 z \in {\bold C}^n, t \in {\bold R}\}.$
We identify $\overline{N_U}$ and $P$ with $\partial
{\Bbb H}_{\bold C}^{n+1}$ and ${\Bbb H}_{\bold C}^{n+1}$
under the map that $\overline{n}(z, t)$ and
$\overline{n}(z, t) a(\zeta)$ are identified with
$(z, t)$ and $(z, t, \rho)$, respectively. Here
$\rho=e^{2 \zeta}$.

  The multiplication of $\overline{N_U}$ (or $\partial
{\Bbb H}_{\bold C}^{n+1}$) is given by
$$(z, t) (z^{\prime}, t^{\prime})=(z+z^{\prime},
  t+t^{\prime}+2 \text{Im} z \overline{z^{\prime}}),$$
where $z \overline{z^{\prime}}=\sum_{j=1}^{n} z_j
\overline{z_j^{\prime}}$. So $\overline{N_U}$ is the
Heisenberg group ${\bold H}^n$. The delation of
${\bold H}^n$ is given by $\rho(z, t)=(\sqrt{\rho} z,
\rho t)$, $\rho>0$, which is consistent with the delation
of ${\Bbb H}_{\bold C}^{n+1}$ given by $\rho(z, z_{n+1})
=(\sqrt{\rho} z, \rho z_{n+1})$. The multiplication is
defined as
$$(z, t, \rho) (z^{\prime}, t^{\prime}, \rho^{\prime})=
  (z+\sqrt{\rho} z^{\prime}, t+\rho t^{\prime}+2
  \sqrt{\rho} \text{Im} z \overline{z^{\prime}},
  \rho \rho^{\prime}).$$
$P_U$ is a locally compact nonunimodular group
with the left Haar measure
$$d \sigma(z, t, \rho)=\rho^{-(n+2)} dm(z) dt d\rho,$$
where $dm(z)$ denotes the Lebesque measure of ${\bold C}^n$.

  From the theory of integrals on quotients $G/H$
where $H$ is a closed subgroup of the Lie group $G$,
we will need an understanding of the formula:
$$\int_{G} f(g) dg=\int_{G/H} \int_{H} f(gh) dh
  d\overline{g}.$$
Here $dg$ and $dh$ are Haar measures on $G$, $H$,
respectively. And this formula defines the $G$-invariant
measure $d\overline{g}$ on the quotient space $G/H$.
Such an integral is determined up to a positive constant.
Formula holds provided that both $G$ and $H$ are unimodular.
So $d \sigma$ is $G$-invariant measure. We have
$$\aligned
 &\text{Vol}(G({\bold Z}) \backslash {\Bbb H}_{\bold C}
  ^{n+1})
 =\int_{G({\bold Z}) \backslash {\Bbb H}_{\bold C}
  ^{n+1}} \rho^{-(n+2)} dt d\rho dm(z)\\
=&\int_{-\frac{1}{2}}^{\frac{1}{2}} dx_{n+1}
  \int_{\sqrt{1-x_{n+1}^2}}^{\infty} dy_{n+1}
  \int_{\sum_{j=1}^{n} (x_j^2+y_j^2)<y_{n+1}}
  \frac{dx_1 dy_1 \cdots dx_n dy_n}{[y_{n+1}-
  \sum_{j=1}^{n} (x_j^2+y_j^2)]^{n+2}}\\
=&\int_{-\frac{1}{2}}^{\frac{1}{2}} dx_{n+1}
  \int_{\sqrt{1-x_{n+1}^2}}^{\infty} dy_{n+1}
  \int_{S^{2n-1}} d\omega
  \int_{0}^{\sqrt{y_{n+1}}}
  \frac{r^{2n-1} dr}{(y_{n+1}-r^2)^{n+2}},
\endaligned$$
where $\int_0^{\sqrt{y_{n+1}}} \frac{r^{2n-1}
dr}{(y_{n+1}-r^2)^{n+2}}=\frac{1}{2} \int_0
^{y_{n+1}} \frac{(y_{n+1}-t)^{n-1}}{t^{n+2}}
dt=\infty$. Therefore, the covolume of
$G({\bold Z})$ is infinite.

\vskip 0.5 cm
\centerline{\bf 3. The eigenfunctions of $L$}
\vskip 0.5 cm

\flushpar
  In this section, we will solve the eigenfunctions of $L$.

  By the transform $u_j=x_j, v_j=y_j (1 \leq j \leq n)$,
$t=x_{n+1}, \rho=y_{n+1}-\sum_{j=1}^n (x_j^2+y_j^2)$ and
Theorem 2.2, we have the following theorem:

\proclaim\nofrills{\smc Theorem 3.1}.  
{\it In coordinates $(x_j, y_j, t, \rho)$, $x_j=\text{Re}z_j$,
$y_j=\text{Im}z_j$, $t=\text{Re}z_{n+1}$, $\rho=\text{Im}z_{n+1}
-\sum_{j=1}^{n} |z_j|^2$, $x_j \ne 0$ and $y_j \ne 0$,
$$\aligned
L=&\rho [\frac{1}{4} \sum_{j=1}^n (\frac{\partial^2}
   {\partial x_j^2}+\frac{\partial^2}{\partial y_j^2})+
   (\rho+\sum_{j=1}^n (x_j^2+y_j^2)) \frac{\partial^2}
   {\partial t^2}+\rho \frac{\partial^2}{\partial \rho^2}
   -n \frac{\partial}{\partial \rho}+\\
  &\sum_{j=1}^n (y_j \frac{\partial}{\partial x_j}
   -x_j \frac{\partial}{\partial y_j}) \frac{\partial}
   {\partial t}].
\endaligned\eqno (3.1)
$$}
\endproclaim

  In particular, if $f=\rho^s$, then
$L f=(\rho^2 \frac{\partial^2}{\partial \rho^2}-n \rho
\frac{\partial}{\partial \rho}) \rho^s=s(s-n-1) \rho^s.$

  The definition of cusp forms on $U(n+1, 1)$ requires
that $\int_0^1 f dt=0$. The corresponding Fourier expansion:
$f=\sum_{a} c(a) Z(a, x_1, y_1, \cdots, x_n, y_n, \rho)
   e^{2 \pi iat}.$
By transform
$x_j=\sqrt{\beta_j} \cos \theta_j$, $y_j=
 \sqrt{\beta_j} \sin \theta_j (1 \leq j \leq n)$,
we have the following theorem:

\proclaim\nofrills{\smc Theorem 3.2}. {\it In coordinates $(\beta_j,
\theta_j, t, \rho)$,
$$L=\rho [\sum_{j=1}^n \beta_j \frac{\partial^2}
    {\partial \beta_j^2}+\rho \frac{\partial^2}
    {\partial \rho^2}+(\rho+\sum_{j=1}^n \beta_j)
    \frac{\partial^2}{\partial t^2}+\sum_{j=1}^n
    \frac{1}{4 \beta_j} \frac{\partial^2}{\partial
    \theta_j^2}-\sum_{j=1}^n \frac{\partial^2}{\partial
    \theta_j \partial t}+\sum_{j=1}^n \frac{\partial}
    {\partial \beta_j}-n \frac{\partial}{\partial
    \rho}]. \eqno (3.2) $$}
\endproclaim

  Set
$$f_{a, b}=u(a, b_1, \cdots, b_n, \rho) v(a, b_1,
\cdots, b_n, \beta_1, \cdots, \beta_n) e^{2 \pi iat}
e^{2i (\sum_{j=1}^{n} b_j \theta_j)},$$
by $Lf_{a, b}=\lambda f_{a, b}$, we have
$$\sum_{j=1}^n \beta_j \frac{1}{v} \frac{\partial^2 v}
  {\partial \beta_j^2}+\rho \frac{u^{\prime \prime}}{u}
  (\rho)-4 \pi^2 a^2 (\rho+\sum_{j=1}^n \beta_j)+
  \sum_{j=1}^n [-\frac{b_j^2}{\beta_j}+4 \pi a b_j+\frac{1}{v} 
  \frac{\partial v}{\partial \beta_j}]-n \frac{u^{\prime}}{u}(\rho)=\frac{\lambda}
  {\rho}.$$
Therefore,
$$\aligned
 &\rho \frac{u^{\prime \prime}}{u}(\rho)-4 \pi^2 a^2
  \rho-n \frac{u^{\prime}}{u}(\rho)-\frac{\lambda}{\rho}
  +4 \pi a \sum_{j=1}^n b_j\\
=&-\sum_{j=1}^n \beta_j \frac{1}{v} \frac{\partial^2
  v}{\partial \beta_j^2}+4 \pi^2 a^2 \sum_{j=1}^n
  \beta_j+\sum_{j=1}^n \frac{b_j^2}{\beta_j}-
  \sum_{j=1}^n \frac{1}{v} \frac{\partial v}
  {\partial \beta_j}=k=\text{const}.
\endaligned
$$
\flushpar
(1) $\rho^2 u^{\prime \prime}(\rho)-n \rho
 u^{\prime}(\rho)-[\lambda+(k-4 \pi a
 \sum_{j=1}^n b_j) \rho+4 \pi^2 a^2 \rho^2] u(\rho)=0.$
 
  Set $k=4 \pi a \sum_{j=1}^n b_j$, then we have
$\rho^2 u^{\prime \prime}(\rho)-n \rho u^{\prime}(\rho)-
 (\lambda+4 \pi^2 a^2 \rho^2) u(\rho)=0.$
Let $u(\rho)=\rho^{\frac{n+1}{2}} w(\rho)$, then
$\rho^2 w^{\prime \prime}(\rho)+\rho w^{\prime}(\rho)-
[\lambda+(\frac{n+1}{2})^2+4 \pi^2 a^2 \rho^2] w(\rho)=0.$
A solution is
$w(\rho)=K_{s-\frac{n+1}{2}}(2 \pi |a| \rho)$, 
$\lambda=s(s-n-1)$.
Here $K$-Bessel function $K_s(z)$ is defined as
$K_s(z)=\frac{1}{2} \int_0^{\infty} \exp[-
  \frac{z}{2}(t+\frac{1}{t})] t^{s-1} dt$,
for $\text{Re}(z)>0$, $\text{Re}(s)>0$.
\flushpar
(2) $\sum_{j=1}^n \beta_j \frac{\partial^2 v}{\partial
\beta_j^2}+\sum_{j=1}^n \frac{\partial v}{\partial
\beta_j}+4 \pi a \sum_{j=1}^n b_j v-\sum_{j=1}^n
\frac{b_j^2}{\beta_j} v-4 \pi^2 a^2 \sum_{j=1}^n
\beta_j v=0$.

  For simplicity, we consider the case that $b_j=0$,
$1 \leq j \leq n$.
\flushpar
(i)$v=v_1(\beta_1) \cdots v_n(\beta_n)$, then
$\sum_{j=1}^n (\beta_j \frac{v_j^{\prime \prime}}{v_j}
(\beta_j)+\frac{v_j^{\prime}}{v_j}(\beta_j)-4 \pi^2 a^2
\beta_j)=0$. A particular case is
$\beta_j v_j^{\prime \prime}(\beta_j)+v_j^{\prime}
(\beta_j)-4 \pi^2 a^2 \beta_j v_j(\beta_j)=0,
1 \leq j \leq n$. 
A solution is
$v_j=K_0(2 \pi |a| \beta_j)$.
\flushpar
(ii) $v=v(\beta)$ with $\beta=\beta_1+\cdots+\beta_n$, then
$\beta v^{\prime \prime}(\beta)+n v^{\prime}(\beta)-
  4 \pi^2 a^2 \beta v(\beta)=0$.
Set $v(\beta)=\beta^{\frac{1-n}{2}} w(\beta)$, then
$\beta^2 w^{\prime \prime}(\beta)+\beta w^{\prime}(\beta)-
 [ (\frac{n-1}{2})^2+4 \pi^2 a^2 \beta^2] w(\beta)=0$.
A solution is
$w(\beta)=K_{\frac{n-1}{2}}(2 \pi |a| \beta)$.

  We have the following theorem:

\proclaim\nofrills{\smc Theorem 3.3}. {\it Two solutions
(we call them the normal solutions) of $Lf=\lambda f$ are
$$\aligned
f &=\sum_{m=0}^{\infty} a_m \rho^{\frac{n+1}{2}} K_{s-
    \frac{n+1}{2}}(2 \pi |m| \rho) \prod_{j=1}^n K_0(2 \pi
    |m| \beta_j) e^{2 \pi imt}, \quad \text{and}\\
g &=\sum_{m=0}^{\infty} b_m \rho^{\frac{n+1}{2}}
    K_{s-\frac{n+1}{2}}(2 \pi |m| \rho)
    \beta^{\frac{1-n}{2}} K_{\frac{n-1}{2}}(2 \pi
    |m| \beta) e^{2 \pi i m t},
\endaligned\eqno (3.3)$$
where $\lambda=s(s-n-1)$.}
\endproclaim

  Set $f=\phi(\beta, t, \rho)$ with $\beta=\sum_{i=1}^{n}
\beta_i$, then
$$L f=\rho [\beta \frac{\partial^2}{\partial
\beta^2}+\rho \frac{\partial^2}{\partial \rho^2}+(\rho+
\beta) \frac{\partial^2}{\partial t^2}+n \frac{\partial}
{\partial \beta}-n \frac{\partial}{\partial \rho}] \phi.$$
Set $\phi_s=\rho^s \beta^{\nu}$, then
$L \phi_s=s(s-n-1) \phi_s+t(t+n-1) \rho^{s+1} \beta^{\nu-1}.$
When $\nu=0$ or $\nu=1-n$, $L \phi_s=s(s-n-1) \phi_s$.
Let $\phi_{s, 0}:=\rho(Z)^s$, $\phi_{s, 1-n}:=\rho(Z)^s
\beta(Z)^{1-n}.$ The Eisenstein series
$$E(Z; s, 0):=\sum_{\gamma \in G({\bold Z})_{\infty}
  \backslash G({\bold Z})} \phi_{s, 0} \circ \gamma,
  \quad E(Z; s, 1-n):=\sum_{\gamma \in G({\bold Z})_{\infty}
  \backslash G({\bold Z})} \phi_{s, 1-n} \circ \gamma
  \eqno (3.4)$$
satisfy the following:

\proclaim\nofrills{\smc Theorem 3.4}. {\it
$$L E(Z; s, 0)=s(s-n-1) E(Z; s, 0), \quad 
  L E(Z; s, 1-n)=s(s-n-1) E(Z; s, 1-n). \eqno (3.5)$$}
\endproclaim

  Let $f=f(\tau, \omega)$ with $\omega=\sqrt{\rho}$ and
$\tau=(x_1, \cdots, x_n, y_1, \cdots, y_n)$, then
$$Lf=
 \frac{1}{4}[\omega^2(\frac{\partial^2}{\partial x_1^2}
 +\cdots+\frac{\partial^2}{\partial x_n^2}+\frac{\partial^2}
 {\partial y_1^2}+\cdots+\frac{\partial^2}{\partial y_n^2}+
 \frac{\partial^2}{\partial \omega^2})-(2n+1) \omega
 \frac{\partial}{\partial \omega}]f. \eqno (3.6)$$
Denote $e(x)=\exp(2 \pi i x)$ and set
$f_a=u(\omega) e(a \cdot \tau)$ with
$a=(a_1, \cdots, a_{2n}) \in {\bold Z}^{2n}$ and
$|a|:=(\sum_{j=1}^{2n} a_j^2)^{\frac{1}{2}}$.
By $L f_a=\lambda f_a$, we have
$\omega^2 u^{\prime \prime}(\omega)-(2n+1) \omega
 u^{\prime}(\omega)-(4 \lambda+4 \pi^2 |a|^2
 \omega^2) u(\omega)=0.$
Let $u(\omega)=\omega^{n+1} w(\omega)$, then
$\omega^2 w^{\prime \prime}(\omega)+\omega
 w^{\prime}(\omega)-[4 \lambda+(n+1)^2+4 \pi^2
 |a|^2 \omega^2] w(\omega)=0.$
A solution is
$w(\omega)=\omega^{n+1} K_{2s-(n+1)}(2 \pi |a| \omega)$
with $\lambda=s(s-n-1)$.

\proclaim\nofrills{\smc Theorem 3.5}. {\it The other solution of
$Lf=\lambda f$ (we call it the singular solution) is
$$f=\sum_{a \in {\bold Z}^{2n}} A(a) \omega^{n+1}
  K_{2s-(n+1)}(2 \pi |a| \omega) e(a \cdot \tau),
  \eqno (3.7)$$
where $\lambda=s(s-n-1)$.}
\endproclaim

  It is known that the Cygan metric $\rho_c$ attached to
$U(n+1, 1)$ is given by (see \cite{1})
$||(z, \rho, t) ||_c=| ||z||^2+\rho-it|^{\frac{1}{2}}$,
for $(z, \rho, t) \in {\bold C}^n \times {\bold R}
\times (0, \infty).$ 
We define the pseudo-distance $d^{*}$ as follows
$$d^{*}(Z, Z^{\prime}):=\log ||(z, t, \rho)^{-1} (z^{\prime},
  t^{\prime}, \rho^{\prime})||_c
=\log \frac{1}{\sqrt{\rho}}| |z-z^{\prime}|^2+\rho^{\prime}
   +i(t-t^{\prime}+2 \text{Im}(z \overline{z^{\prime}}))|^
   {\frac{1}{2}}.$$

\proclaim\nofrills{\smc Theorem 3.6}. {\it The distance function is
noneuclidean harmonic, i.e.,
$$L d^{*}((z^{\prime}, t^{\prime}, \rho^{\prime}),
(z, t, \rho))=0. \eqno (3.8) $$}
\endproclaim

\vskip 0.5 cm
\centerline{\bf 4. The integral transform of
Eisenstein series on $U(n+1, 1)$}
\vskip 0.5 cm

\flushpar
  In this section we will prove the following theorem:

\proclaim\nofrills{\smc Theorem 4.1}. {\it (The integral transform property
of Eisenstein series on $U(n+1, 1)$): The Eisenstein series
have the following Fourier expansions:
$$E(Z; s, 0)=\sum_{m} a_m(\rho, \beta) e^{2 \pi imt}, \quad
  E(Z; s, 1-n)=\sum_{m} b_m(\rho, \beta) e^{2 \pi imt}.$$
(1)Let
$$a_m(\rho):=\int_{({\bold R}_{+})^n} [a_m(\rho, \beta)
  -\delta_{0, m} \rho^s] \prod_{j=1}^n e^{-2 \pi |m| \beta_j}
  G(2 \pi |m| \beta_j) d\beta_1 \cdots d\beta_n,$$
where 
$G(z)=\sum_{k=0}^{\infty} \frac{1}{2^k k!} C_{2k}^k z^k$
and $\delta_{i, j}$ is the Kronecker symbol. Then 
$$a_m(\rho)=\left\{\aligned
 &2^{1-n} \pi^{s-\frac{n}{2}} \frac{\Gamma(s-\frac{n}{2})}
  {\Gamma(s)^2} \varphi_m(s) |m|^{s-\frac{n+1}{2}}
  \rho^{\frac{n+1}{2}} K_{s-\frac{n+1}{2}}
  (2 \pi |m| \rho), (m \ne 0),\\
 &2^{-n} \sqrt{\pi} \frac{\Gamma(s-\frac{n}{2}) \Gamma(s-
  \frac{n+1}{2})}{\Gamma(s)^2} \varphi_0(s) \rho^{n+1-s},
  (m=0),
  \endaligned\right. \eqno (4.1)$$
with $\varphi_m(s)=\sum_{c>0} \frac{1}{|c|^{2s}}
(\sum_{(d, c)=1, d \text{mod} c} e(\frac{md}{c}))$,
where $e(x)$ denotes $\exp(2 \pi i x)$.
\flushpar
(2)For $m \ne 0$, let
$$b_m(\rho)=\int_0^{\infty} b_m(\rho, \beta) e^{-2 \pi
 |m| \beta} H(2 \pi |m| \beta) d\beta,$$
where $H(z)=\sum_{k=n-1}^{\infty} \frac{2^k (1-\frac{n}{2},
k)} {k! (k-n+1)!} z^k,$ 
then 
$$\aligned
 b_m(\rho)=
&2^{n-1} \pi^{s-\frac{n}{2}}
  \frac{\Gamma(s-\frac{n}{2})}{\Gamma(s-n+1)^2}
  \varphi_m(s-n+1) |m|^{s-\frac{n+1}{2}}
  \rho^{\frac{n+1}{2}} K_{s-\frac{n+1}{2}}(2
  \pi |m| \rho)\\
&-R(|m|, \rho)
\endaligned\eqno (4.2)$$
with
$$\aligned
R(|m|, \rho)= &2^{n-2} \pi^{s-\frac{n}{2}}
\Gamma(s-n+1)^{-1} \varphi_m(s-n+1) |m|^{s-1-\frac{n}{2}}
\rho^{\frac{n}{2}}\\
&\times \sum_{k=0}^{n-2} \frac{(1-\frac{n}{2}, k)}{k!}
(4 \pi |m| \rho)^{\frac{k}{2}} W_{\frac{n}{2}-1-\frac{k}{2},
s-\frac{n+1}{2}-\frac{k}{2}}(4 \pi |m| \rho).
\endaligned\eqno (4.3)
$$
Here $(\alpha, k):=\alpha (\alpha+1) \cdots (\alpha+k-1)$
and $W_{\kappa, \mu}(x)$ is the Whittaker function.}
\endproclaim

\demo\nofrills{Proof}. We have
$$
\aligned
   a_m(\rho, \beta)
= &\int_0^1 \frac{1}{2} \sum_{c, d \in {\bold Z}, (c, d)=1}
   \frac{\rho^s}{[(ct+d)^2+c^2 (\rho+\beta)^2]^s}
   e^{-2 \pi imt} dt\\
= &\delta_{0, m} \rho^s+\sum_{c>0} \frac{1}{|c|^{2s}}
   \int_{-\infty}^{\infty} \sum_{(d, c)=1, d \text{mod} c}
   \frac{\rho^s}{[(t+\frac{d}{c})^2+(\rho+\beta)^2]^s}
   e^{-2 \pi imt} dt\\
= &\delta_{0, m} \rho^s+\sum_{c>0} \frac{1}{|c|^{2s}}
   \sum_{(d, c)=1, d \text{mod} c} e(\frac{md}{c})
   \int_{-\infty}^{\infty} \frac{\rho^s e^{-2 \pi imt}}
   {[t^2+(\rho+\beta)^2]^s} dt.
\endaligned
$$

  Set $\varphi_m(s)=\sum_{c>0} \frac{1}{|c|^{2s}}
(\sum_{(d, c)=1,  d\text{mod} c} e(\frac{md}{c}))$,
by (see \cite{7}, p. 15)
$$\int_{-\infty}^{\infty} \frac{e(-ut)}{(1+t^2)^s} dt=
   2 \pi^s |u|^{s-\frac{1}{2}} \Gamma(s)^{-1} 
   K_{s-\frac{1}{2}}(2 \pi |u|), \quad
   (u \ne 0, u \in {\bold R}), \eqno (4.4)$$
we have for $m \ne 0$,
$$\int_{-\infty}^{\infty} \frac{\rho^s e^{-2 \pi imt}}
  {[t^2+(\rho+\beta)^2]^s} dt=\frac{2 \pi^s}{\Gamma(s)}
  |m|^{s-\frac{1}{2}} \rho^s (\rho+\beta)^{\frac{1}{2}-s}
  K_{s-\frac{1}{2}}(2 \pi |m| (\rho+\beta)).$$
So,
$$
\aligned
a_m(\rho)
 =&\varphi_m(s) \frac{2 \pi^s}{\Gamma(s)} |m|^{s-\frac{1}{2}}
  \rho^s \int_0^{\infty} \cdots \int_0^{\infty}
  (\rho+\beta)^{\frac{1}{2}-s}\\
  &K_{s-\frac{1}{2}}(2 \pi |m| (\rho+\beta))
  \prod_{j=1}^n e^{-2 \pi |m| \beta_j} G(2 \pi |m| \beta_j)
  d\beta_1 \cdots d\beta_n.
\endaligned
$$
Denote the above integral as $A_m(\rho)$, set
$\rho_j=\rho+\beta_{j+1}+\cdots+\beta_n$ for $1 \leq j \leq n$,
then
$$\aligned
  A_m(\rho)
=&\int_0^{\infty} \cdots \int_0^{\infty} e^{-2 \pi |m|
  (\beta_2+\cdots+\beta_n)} G(2 \pi |m| \beta_2) \cdots
  G(2 \pi |m| \beta_n)\\
 &\int_0^{\infty} (\rho_1+\beta_1)^{\frac{1}{2}-s}
  K_{s-\frac{1}{2}}(2 \pi |m| (\rho_1+\beta_1))
  e^{-2 \pi |m| \beta_1} G(2 \pi |m| \beta_1) d \beta_1
  \cdots d \beta_n.
\endaligned
$$

  We recall that the Weyl fractional integral is defined as
$$h(y; \mu)=\frac{1}{\Gamma(\mu)} \int_y^{\infty} f(x)
  (x-y)^{\mu-1} dx. \eqno (4.5)$$
When $f(x)=x^{-\nu} e^{-\alpha x} K_{\nu}(\alpha x)$ and
$\text{Re}(\mu)>0$,
$$h(y; \mu)=\sqrt{\pi} (2 \alpha)^{-\frac{1}{2} \mu-
  \frac{1}{2}} y^{\frac{1}{2} \mu-\nu-\frac{1}{2}}
  e^{-\alpha y} W_{-\frac{1}{2} \mu, \nu-\frac{1}{2}
  \mu}(2 \alpha y), \eqno (4.6)$$
for $\text{Re}(\alpha y)>0$ (see \cite{3}, p. 208, (53)).
Here the Whittaker functions (see \cite{4}, Vol. I, p. 264)
$W_{\kappa, \mu}(x)=e^{-\frac{x}{2}} x^{\frac{c}{2}} \Psi(a, c; x),$
where $a=\frac{1}{2}-\kappa+\mu$, $c=2 \mu+1$ and
$\Psi(a, c; x)=\frac{1}{\Gamma(a)} \int_0^{\infty}
  e^{-xt} t^{a-1} (1+t)^{c-a-1} dt$,
$\text{Re}(a)>0$. Thus
$$W_{\kappa, \mu}(x)=e^{-\frac{x}{2}} x^{\mu+\frac{1}{2}}
  \frac{1}{\Gamma(\frac{1}{2}-\kappa+\mu)} \int_0^{\infty}
  e^{-xt} t^{-\frac{1}{2}-\kappa+\mu} (1+t)^{-\frac{1}{2}+
  \kappa+\mu} dt.$$

  It is known that (see \cite{4}, Vol. I, p. 265, (13)) 
$$\aligned
 K_{\nu}(x)
&=\sqrt{\pi} e^{-x} (2x)^{\nu} \Psi(\frac{1}{2}
  +\nu, 1+2 \nu; 2x)\\
&=\frac{\sqrt{\pi} e^{-x} (2x)^{\nu}}
  {\Gamma(\nu+\frac{1}{2})}\int_0^{\infty} e^{-2xt}
  t^{\nu-\frac{1}{2}}(1+t)^{\nu-\frac{1}{2}} dt.
\endaligned
$$

  By the above formulas, we have
$$\aligned
   &A_{1, m}(\rho_1) \\
:= &\int_0^{\infty} (\rho_1+\beta_1)^{\frac{1}{2}-s}
  K_{s-\frac{1}{2}}(2 \pi |m| (\rho_1+\beta_1))
  e^{-2 \pi |m| \beta_1} G(2 \pi |m| \beta_1) d \beta_1\\
=&e^{2 \pi |m| \rho_1} \sum_{k=0}^{\infty} \frac{1}{2^k k!}
  C_{2k}^k (2 \pi |m|)^k \int_{\rho_1}^{\infty}
  \beta_1^{\frac{1}{2}-s} K_{s-\frac{1}{2}}(2 \pi |m| \beta_1)
  e^{-2 \pi |m| \beta_1} (\beta_1-\rho_1)^k d\beta_1\\
=&\frac{\sqrt{\pi}}{4 \pi |m|} \rho_1^{\frac{1}{2}-s}
  \sum_{k=0}^{\infty} \frac{1}{2^k} C_{2k}^k (\pi |m|
  \rho_1)^{\frac{k}{2}} W_{-\frac{k}{2}-\frac{1}{2},
  s-1-\frac{k}{2}}(4 \pi |m| \rho_1)\\
=&\frac{\sqrt{\pi}}{4 \pi |m|} \frac{1}{\Gamma(s)}
  (4 \pi |m|)^{s-\frac{1}{2}} e^{-2 \pi |m| \rho_1}
  \int_0^{\infty} e^{-4 \pi |m| \rho_1 t} (1+t)^{s-2} t^{s-1}\\
 &\times \sum_{k=0}^{\infty} \frac{1}{2^{2k}} C_{2k}^k (1+t)^{-k} dt\\
=&\frac{\sqrt{\pi}}{4 \pi |m| \Gamma(s)}
  (4 \pi |m|)^{s-\frac{1}{2}} e^{-2 \pi |m| \rho_1}
  \int_0^{\infty} e^{-4 \pi |m| \rho_1 t}
  (1+t)^{s-\frac{3}{2}} t^{s-\frac{3}{2}} dt\\
=&\frac{\Gamma(s-\frac{1}{2})}{\sqrt{4 \pi |m|} \Gamma(s)}
  \rho_1^{1-s} K_{s-1}(2 \pi |m| \rho_1).
\endaligned
$$

  In general, we have
$$\aligned
  A_{j, m}(\rho_j) :=  
 &\int_0^{\infty} (\rho_j+\beta_j)^{\frac{j}{2}-s}
  K_{s-\frac{j}{2}}(2 \pi |m| (\rho_j+\beta_j))
  e^{-2 \pi |m| \beta_j} G(2 \pi |m| \beta_j) d\beta_j\\
=&\frac{\Gamma(s-\frac{j}{2})}{\sqrt{4 \pi |m|}
  \Gamma(s-\frac{j}{2}+\frac{1}{2})} \rho_j^{\frac{j+1}
  {2}-s} K_{s-\frac{j+1}{2}}(2 \pi |m| \rho_j),
\endaligned
$$
for $1 \leq j \leq n$ and $\rho_n=\rho$. Therefore,
$$\aligned
  A_m(\rho)
 =&\prod_{j=1}^{n} \frac{\Gamma(s-\frac{j}{2})}{\sqrt{4 \pi |m|}
   \Gamma(s-\frac{j}{2}+\frac{1}{2})} \rho^{\frac{n+1}{2}-s}
   K_{s-\frac{n+1}{2}}(2 \pi |m| \rho)\\
 =&(4 \pi |m|)^{-\frac{n}{2}} \frac{\Gamma(s-
   \frac{n}{2})}{\Gamma(s)} \rho^{\frac{n+1}{2}-s}
   K_{s-\frac{n+1}{2}}(2 \pi |m| \rho).
\endaligned
$$
Similarly, the other part of the theorem can be proved.
\quad \qquad \qquad \qquad \qquad \qquad $\boxed{}$
\enddemo

\vskip 0.5 cm
\centerline{\bf 5. The Poincar\'{e} series for $U(n+1, 1)$}
\vskip 0.5 cm

\flushpar
  The concept of a point-pair invariant was introduced
by Selberg \cite{9} who made fascinating use of it. Now, we
introduce the following concept.

\definition{\it Definition 5.1} A map $f: {\Bbb H}_{\bold C}^{n+1}
\times {\Bbb H}_{\bold C}^{n+1} \to {\bold C}$ is called {\it a
point-pair invariant associated to a discrete subgroup \
$\Gamma \leq G({\bold Z})$} if
$f(\gamma(P), \gamma(Q))=f(P, Q)$ for all $P, Q \in
{\Bbb H}_{\bold C}^{n+1}$ and $\gamma \in \Gamma$.
\enddefinition

\proclaim\nofrills{\smc Theorem 5.2}.
{\it For $L=\rho [\beta \frac{\partial^2}{\partial \beta^2}+
\rho \frac{\partial^2}{\partial \rho^2}+(\rho+\beta)
\frac{\partial^2}{\partial t^2}+n \frac{\partial}{\partial
\beta}-n \frac{\partial}{\partial \rho}]$, set
$f=g(u, v)$, where
$$u=u(Z, Z^{\prime}):=\frac{(t-t^{\prime})^2+(\rho+\beta-
  \rho^{\prime}-\beta^{\prime})^2}{4 \rho \rho^{\prime}},
  \quad v=v(Z, Z^{\prime}):=\frac{\beta \beta^{\prime}}{\rho
  \rho^{\prime}} \eqno (5.1)$$
are two point-pair invariants associated to $G({\bold Z})$.
Then
$$\aligned
L f=&[(u^2+(\lambda+1)u) \frac{\partial^2}{\partial u^2}+
     2uv \frac{\partial^2}{\partial u \partial v}
     +v(v+\lambda) \frac{\partial^2}{\partial v^2}\\
    &+((n+2)u+\lambda+1) \frac{\partial}{\partial u}
     +((n+2)v+n \lambda) \frac{\partial}{\partial v}] g
\endaligned\eqno (5.2)$$
with $\lambda=\frac{\beta^{\prime}}{\rho^{\prime}}$,
where $\gamma(\lambda)=\lambda$ for $\gamma \in
G({\bold Z})$.}
\endproclaim

\demo\nofrills{Proof}. It is obtained by a straightforward calculation.
\qquad \qquad \qquad \qquad \qquad \qquad $\boxed{}$
\enddemo

  Let
$M(u, v, \frac{\partial}{\partial u}, \frac{\partial}
{\partial v})=[u^2+(\lambda+1)u] \frac{\partial^2}{\partial
u^2}+2uv \frac{\partial^2}{\partial u \partial v}
+v(v+\lambda) \frac{\partial^2}{\partial v^2}+[(n+2)u+\lambda+1]
\frac{\partial}{\partial u}+[(n+2)v+n \lambda ] \frac{\partial}
{\partial v}+s(n+1-s).$
By transform $x=\frac{u}{\lambda+1}, y=\frac{v}{\lambda}$, we
have
$$\aligned
  M(x, y, \frac{\partial}{\partial x},
  \frac{\partial}{\partial y})
=&x(x+1) \frac{\partial^2}{\partial x^2}+2xy \frac{\partial^2}
{\partial x \partial y}+y(y+1) \frac{\partial^2}{\partial y^2}\\
&+[(n+2)x+1]\frac{\partial}{\partial x}+[(n+2)y+n]
\frac{\partial}{\partial y}+s(n+1-s).
\endaligned\eqno (5.3)$$

\proclaim\nofrills{\smc Theorem 5.3}. {\it Some solutions of the equation
$M(x, y, \frac{\partial}{\partial x}, \frac{\partial}
{\partial y}) g(x, y)=0$ are as follows:
$$\aligned
 &g_s(x, y)=x^{-a} y^{-b} F_3(a, b; a, b-n+1; 2s-n;
  -x^{-1}, -y^{-1}),\\
 &g_1(x)=x^{-s} {}_{2}F_{1}(s, s; 2s-n; -x^{-1}),\\
 &g_2(y)=y^{-s} {}_{2}F_{1}(s, s-n+1; 2s-n; -y^{-1}),\\
 &g_3(x+y)=w^{-s} {}_{2}F_{1}(s, s-n; 2s-n; -(x+y)^{-1}).
\endaligned\eqno (5.4)$$}
\endproclaim

\demo\nofrills{Proof}.
At first, we consider the degenerate case.
\flushpar
(1) Set $g(x, y)=g_1(x)$, then one has
$$M(x, \frac{d}{d x}) g_1(x)=\{ x(x+1)
  \frac{d^2}{d x^2}+[(n+2)x+1] \frac{d}{d x}
  +s(n+1-s) \} g_1(x)=0.$$
A solution is
$g_1(x)=x^{-s} {}_{2}F_{1}(s, s; 2s-n; -x^{-1})$.
\flushpar
(2) Set $g(x, y)=g_2(y)$, then one has
$$M(y, \frac{d}{d y}) g_2(y)=\{ y(y+1)
  \frac{d^2}{d y^2}+[(n+2)y+n] \frac{d}{d y}
  +s(n+1-s) \} g_2(y)=0.$$
A solution is
$g_2(y)=y^{-s} {}_{2}F_{1}(s, s-n+1; 2s-n; -y^{-1})$.
\flushpar
(3) Set $g(x, y)=g_3(w)$ with $w=x+y$, then one has
$$M(w, \frac{d}{d w}) g_3(w)=\{ w(w+1)
  \frac{d^2}{d w^2}+[(n+2)w+(n+1)] \frac{d}{d w}
  +s(n+1-s) \} g_3(w)=0.$$
A solution is
$g_3(w)=w^{-s} {}_{2}F_{1}(s, s-n; 2s-n; -w^{-1})$.

  Secondly, we consider the general case.
Set $g(x, y)=x^{-a} y^{-b} f(-x^{-1}, -y^{-1})$, then
$f(x, y)$ satisfies the following equation:
$$\aligned
 &x \{ x(1-x) \frac{\partial^2 f}{\partial x^2}+
  y \frac{\partial^2 f}{\partial x \partial y}+
  [2(a+b)-n-(2a+1)x] \frac{\partial f}{\partial x}-
  a^2 f \}+\\
 &y \{ y(1-y) \frac{\partial^2 f}{\partial y^2}+
  x \frac{\partial^2 f}{\partial x \partial y}+
  [2(a+b)-n-(2b-n+2)y] \frac{\partial f}{\partial y}-
  b(b-n+1) f \}+\\
 &(a+b-s)[a+b-(n+1-s)] f=0.
\endaligned$$
Let
$$\left\{\aligned
 &x(1-x) \frac{\partial^2 f}{\partial x^2}+
  y \frac{\partial^2 f}{\partial x \partial y}+
  [2(a+b)-n-(2a+1)x] \frac{\partial f}{\partial x}-
  a^2 f=0,\\
 &y(1-y) \frac{\partial^2 f}{\partial y^2}+
  x \frac{\partial^2 f}{\partial x \partial y}+
  [2(a+b)-n-(2b-n+2)y] \frac{\partial f}{\partial y}-
  b(b-n+1) f=0,\\
 &(a+b-s)[a+b-(n+1-s)]=0.
\endaligned\right.$$
Without loss of generality, we can assume that $a+b=s$.
Let $z=f$, $p=\frac{\partial z}{\partial x}, q=\frac{
\partial z}{\partial y}, r=\frac{\partial^2 z}{\partial x^2},
s=\frac{\partial^2 z}{\partial x \partial y},
t=\frac{\partial^2 z}{\partial y^2}$, we have
$$\left\{\aligned
 &x(1-x)r+ys+[2(a+b)-n-(2a+1)x]p-a^2 z=0,\\
 &y(1-y)t+xs+[2(a+b)-n-(2b-n+2)y]q-b(b-n+1) z=0.
\endaligned\right.$$

  It is known that a solution of the equations
$$\left\{\aligned
 &x(1-x)r+ys+[\gamma-(\alpha+\beta+1)x]p-\alpha \beta z=0,\\
 &y(1-y)t+xs+[\gamma-(\alpha^{\prime}+\beta^{\prime}+1)y]q
  -\alpha^{\prime} \beta^{\prime} z=0.
\endaligned\right.$$
is $z=F_3(\alpha, \alpha^{\prime}; \beta, \beta^{\prime};
\gamma; x, y)$(see \cite{2}).
Here
$$\left\{\aligned
 &\gamma=2(a+b)-n,\\
 &\alpha+\beta+1=2a+1, \quad \alpha \beta=a^2,\\
 &\alpha^{\prime}+\beta^{\prime}+1=2b-n+2, \quad
  \alpha^{\prime} \beta^{\prime}=b(b-n+1).
\endaligned\right.$$
i.e., $\alpha=\beta=a$, $\alpha^{\prime}=b$,
$\beta^{\prime}=b-n+1$, $\gamma=2s-n$. Therefore,
a family of solutions are
$$g_s(x, y)=x^{-a} y^{-b} F_3(a, b; a, b-n+1; 2s-n;
  -x^{-1}, -y^{-1}).$$
This completes the proof of Theorem 5.3.
\qquad \quad \qquad \qquad \qquad \qquad \qquad \qquad \qquad $\boxed{}$
\enddemo

\definition{\it Definition 5.4} A function $f: {\Bbb H}_{\bold C}^{n+1}
\to {\bold C}$ is called {\it a nonholomorphic automorphic form
attached to the unitary group $U(n+1, 1)$} if it satisfies
the following three conditions:
\roster
\item
$f$ is an eigenfunction of the Laplace-Beltrami operator
of $U(n+1, 1)$ on ${\Bbb H}_{\bold C}^{n+1}$;
\item
$f$ is invariant under the modualr group; i.e.,
$f(\gamma(Z))=f(Z)$ for all
$\gamma \in G({\bold Z})$ and all $Z \in {\Bbb H}_{\bold C}^{n+1}$.
\item
$f$ has at most polynomial growth at infinity; i.e.,
there are constants $C>0$ and $k$ such that
$|f(Z)| \leq C \rho^k$, as $\rho \to \infty$
uniformly in $t$, for fixed $\beta$.
\endroster
\enddefinition

  We denote by ${\Cal N}(G({\bold Z}), \lambda)$ the
space of such nonholomorphic automorphic forms attached
to $U(n+1, 1)$.

  By Theorem 3.4, we have $E(Z; s, 0), E(Z; s, 1-n) \in 
{\Cal N}(G({\bold Z}), s(s-n-1))$ when $\text{Re}(s)>n+1$.

  Now, we study the structure of
${\Cal N}(G({\bold Z}), s(s-n-1))$.
The Poincar\'{e} series is defined as
$$r(Z, Z^{\prime}; s):=\sum_{\gamma \in G({\bold Z})}
  g_s(x(Z, \gamma(Z^{\prime})), y(Z)). \eqno (5.5)$$

\proclaim\nofrills{\smc Theorem 5.5}. {\it
$r(Z, Z^{\prime}; s) \in {\Cal N}(G({\bold Z}),
s(s-n-1))$ for $\text{Re}(s)>n$, $\text{Re}(a)>1$ and
$\text{Re}(b)>n-1$, where $a+b=s$.}
\endproclaim

\demo\nofrills{Proof}. Without loss of generality, we can only
consider the nondegenerate case.

  According to \cite{2}, the two variable hypergeometric
function $F_3$ has the following integral representation:
$$\aligned
  F_3(\alpha, \alpha^{\prime}; \beta, \beta^{\prime};
  \gamma; x, y)=
 &\frac{\Gamma(\gamma)}{\Gamma(\beta) \Gamma(\beta^{\prime})
  \Gamma(\gamma-\beta-\beta^{\prime})}
  \iint_{u \geq 0, v \geq 0, u+v \leq 1} u^{\beta-1}\\
 &v^{\beta^{\prime}-1} (1-u-v)^{\gamma-\beta-
  \beta^{\prime}-1} (1-ux)^{-\alpha} (1-vy)^
  {-\alpha^{\prime}} du dv,
\endaligned$$
for $\text{Re}(\beta)>0$, $\text{Re}(\beta^{\prime})>0$
and $\text{Re}(\gamma-\beta-\beta^{\prime})>0$.

  Now, one has
$$\aligned
  g_s(x, y)=
 &\frac{\Gamma(2s-n)}{\Gamma(a) \Gamma(b-n+1) \Gamma(s-1)}
  \iint_{u \geq 0, v \geq 0, u+v \leq 1} u^{a-1}\\
 &v^{b-n} (1-u-v)^{s-2} (x+u)^{-a} (y+v)^{-b} du dv,
\endaligned\eqno (5.6)$$
for $\text{Re}(a)>0$, $\text{Re}(b)>n-1$ and $\text{Re}(s)>1$.
On the other hand,
$x=x(Z, Z^{\prime})=(1+\frac{\beta}{\rho}) \sigma(z_{n+1},
z_{n+1}^{\prime})$, where $\sigma(z_{n+1}, z_{n+1}^{\prime})
=\frac{|z_{n+1}-z_{n+1}^{\prime}|^2}{4 \text{Im} z_{n+1}
\text{Im} z_{n+1}^{\prime}}$, and $y=y(Z)=\frac{\beta}{\rho}$.
Therefore, we have
$$\aligned
  r(Z, Z^{\prime}; s)
=&\frac{\Gamma(2s-n)}{\Gamma(a) \Gamma(b-n+1)
  \Gamma(s-1)} \iint_{u \geq 0, v \geq 0, u+v \leq 1}
  u^{a-1} v^{b-n} (y+v)^{-b}\\
 &\times (1-u-v)^{s-2} \sum_{\gamma \in G({\bold Z})}
  \frac{1}{[(1+\frac{\beta}{\rho}) \sigma(z_{n+1},
  \gamma(z_{n+1}^{\prime}))+u]^{a}} du dv.
\endaligned\eqno (5.7)$$
The sum
$$\sum_{\gamma \in G({\bold Z})} \frac{1}{[(1+\frac
  {\beta}{\rho}) \sigma(z_{n+1}, \gamma(z_{n+1}^
  {\prime}))+u]^{\text{Re}(a)}}
\leq \sum_{\gamma \in G({\bold Z})} \frac{1}{[\sigma(
  z_{n+1}, \gamma(z_{n+1}^{\prime}))+u]^{\text{Re}(a)}}.$$
By \cite{8}, p. 285, Lemma 1, if $\text{Re}(a)>1$, then the
series
$$\sum_{\gamma \in G({\bold Z})} \frac{1}{[1+\sigma(
  z_{n+1}, \gamma(z_{n+1}^{\prime}))]^{\text{Re}(a)}}$$
is convergent uniformly for $z_{n+1}$, $z_{n+1}^{\prime}$
in compact domains. By the same method in \cite{8}, one has
that if $z_{n+1} \notin G({\bold Z}) z_{n+1}^{\prime}$,
then the series in (5.7) is convergent absolutely for
$\text{Re}(a)>1$, i.e., $r(Z, Z^{\prime}; s)$
is well-defined for $\text{Re}(a)>1$.

  For $g \in G({\bold Z})$, by $x(\gamma(Z), \gamma(Z^{
\prime}))=x(Z, Z^{\prime})$ and $y(\gamma(Z))=y(Z)$, we have
$$\aligned
 &r(g(Z), Z^{\prime}; s)
=\sum_{\gamma \in G({\bold Z})} g_s(x(g(Z), \gamma(
  Z^{\prime})), y(g(Z)))\\
=&\sum_{\gamma \in G({\bold Z})} g_s(x(g(Z), g \circ \gamma(Z
  ^{\prime})), y(Z))
=\sum_{\gamma \in G({\bold Z})} g_s(x(Z, \gamma(Z^{\prime})),
  y(Z))\\
=&r(Z, Z^{\prime}; s).
\endaligned$$

  For $\gamma \in G({\bold Z})$,
$M(x, y, \frac{\partial}{\partial x}, \frac{\partial}
 {\partial y})(Z, \gamma(Z^{\prime})) g_s(x(Z, \gamma(Z^{
 \prime})), y(Z))=0.$
Hence,
$$(L-s(s-n-1)) g_s(x(Z, \gamma(Z^{\prime})), y(Z))=0.$$
Thus,
$(L-s(s-n-1)) r(Z, Z^{\prime}; s)=0$.
\quad \qquad \qquad \qquad \qquad \qquad \qquad \qquad \qquad \qquad $\boxed{}$
\enddemo

  Theorem 5.5 implies the following theorem:

\proclaim\nofrills{\smc Theorem 5.6}. {\it There exist infinitely many
elements in ${\Cal N}(G({\bold Z}), s(s-n-1))$.}
\endproclaim

\vskip 0.5 cm
\centerline{\bf 6. The Poisson kernel and Eisenstein
series for $U(n+1, 1)$}
\vskip 0.5 cm

\flushpar
  In this section, we will give the Poisson kernel of
$L$ on ${\Bbb H}_{\bold C}^{n+1}$ and the corresponding
Eisenstein series.

  Let us give the Iwasawa decomposition of $G$.
$G=KAN=\overline{N}AK$, where
$$A=T A_U T^{-1}=\left\{
  a=\left(\matrix
  I_n &           &           \\   
      & e^{\zeta} &           \\
      &           & e^{-\zeta}
\endmatrix\right): \zeta \in {\bold R} \right\},$$
$$N=T N_U T^{-1}=\left\{
  n=\left(\matrix
  I_n             & 0 & z^{t}   \\
  2i \overline{z} & 1 & t+i|z|^2\\
  0               & 0 & 1 
\endmatrix\right): t \in {\bold R}, z \in {\bold C}^n
\right\},$$
$$\overline{N}=T \overline{N_U} T^{-1}=
\left\{\overline{n}=\left(\matrix
  I_n           & i z^{t}   & 0\\
  0             & 1          & 0\\
 -2\overline{z} & -t-i |z|^2 & 1
\endmatrix\right): t \in {\bold R}, z \in {\bold C}^n
\right\},$$
and $K=T K_U T^{-1}$.

  For
$\overline{n}=\left(\matrix
  I_n           & i z^{t}   & 0\\
  0             & 1          & 0\\
 -2\overline{z} & -t-i |z|^2 & 1
\endmatrix\right)$,
$\overline{n}^{\prime}
=\left(\matrix
  I_n           & i w^{t}   & 0\\
  0             & 1          & 0\\
 -2\overline{w} & -\zeta-i |w|^2 & 1
\endmatrix\right)$ and
$a=\left(\matrix
  I_n &                     & \\
      & \rho^{-\frac{1}{2}} & \\
      &  & \rho^{\frac{1}{2}}
\endmatrix\right),$
$$[(\overline{n}a)^{-1}]^{*}=\left(\matrix
  I_n & 0 & 2\rho^{-\frac{1}{2}} z^{t}\\
  i \overline{z} & \rho^{\frac{1}{2}} &
  \rho^{-\frac{1}{2}}(t+i|z|^2)\\
  0 & 0 & \rho^{-\frac{1}{2}}
\endmatrix\right).$$
For $Z \in {\Bbb H}_{\bold C}^{n+1}$, consider
$Z^{\prime}=(z^{\prime}, z_{n+1}^{\prime})=(-\frac{z}
{\sqrt{\rho}}, i \frac{y_{n+1}}{\rho})$,
$\text{Im}z_{n+1}^{\prime}-|z^{\prime}|^2=1$, so
$Z^{\prime} \in {\Bbb H}_{\bold C}^{n+1}$. Thus
$$[(\overline{n}a)^{-1}]^{*}(-\frac{z}{\sqrt{\rho}},
  i \frac{y_{n+1}}{\rho})=(z, t+i y_{n+1})=Z.
  \eqno (6.1)$$
$$a^{-1} \overline{n}^{-1} \overline{n}^{\prime}=
 \left(\matrix
  I_n         & -i(z^t-w^t)        & 0\\
  0           & \rho^{\frac{1}{2}} & 0\\
 2\rho^{-\frac{1}{2}}(\overline{z}-\overline{w})
 &\rho^{-\frac{1}{2}}[t-\zeta-i(|z|^2+|w|^2-2
  \overline{z}w^{t})] & \rho^{-\frac{1}{2}}
\endmatrix\right).\eqno (6.2)$$ 

  For $a \in A$ and $n \in N$ and
$k=T k_U T^{-1}=\frac{1}{2i}\left(\matrix
    A_1 & *\\
    A_2 & *
  \endmatrix\right),$
where $k_U=\left(\matrix
          A & \\
            & D
\endmatrix\right)$, $A=(a_{ij}) \in U(n+1)$,
$D \in U(1)$, $A_1$ is an $(n+1) \times (n+1)$-matrix
and $A_2$ is an $1 \times (n+1)$-matrix.
$an=\left(\matrix
     B & *\\
     0 & *
    \endmatrix\right)$ with 
$B=\left(\matrix
     I_n & 0\\
     2ie^{\zeta} \overline{z} & e^{\zeta}
 \endmatrix\right)$.
$kan=\frac{1}{2i} \left(\matrix
       A_1 B & *\\
       A_2 B & *
     \endmatrix\right)$, where
$A_1 B=\left(\matrix
     * & *\\
     * & i e^{\zeta} (a_{n+1, n+1}+D)
\endmatrix\right)$ and
$A_2 B=(*, (-a_{n+1, n+1}+D) e^{\zeta})$.

  For $g=(g_{ij}) \in G$, we have
$$g_{n+1, n+1}=\frac{1}{2} e^{\zeta} (a_{n+1, n+1}+D),
\quad
g_{n+2, n+1}=\frac{1}{2i} e^{\zeta}(-a_{n+1, n+1}+D).$$
$D \in U(1)$ implies that
$e^{2 \zeta}=|g_{n+1, n+1}+i g_{n+2, n+1}|^2$.
Therefore, the Poisson kernel
$$\aligned
  P(Z, W)
 &=|\rho^{\frac{1}{2}}+i \rho^{-\frac{1}{2}}
   [t-\zeta-i(|z|^2+|w|^2-2 \overline{z} w^{t})]|^{-2}\\
 &=\frac{\rho}{|\rho+|z-w|^2+i(t-\zeta-2 \text{Im}
   \overline{z} w^t)|^2},
\endaligned\eqno (6.3)$$
where $W=(w_1, \cdots, w_n, \zeta, 0) \in {\bold C}^n
\times {\bold R}$.
We have $L P(Z, W)^s=s(s-n-1) P(Z, W)^s$. In fact, by
Helgason's conjecture, which was proved by Kashiwara et
al. in \cite{6}, that the eigenfunctions on Riemannian
symmetric spaces can be represented as Poisson integrals
of their hyperfunction boundary values.

  In the theory of automorphic forms, the rigid
property is essential, it is determined by the
discrete subgroup. $G({\bold Z})$ acts on
${\Bbb H}_{\bold C}^{n+1}$, not ${\bold C}^n
\times {\bold R} \times (0, \infty)$, although they
are diffeomorphic. The boundary of
${\Bbb H}_{\bold C}^{n+1}$ is
$$\partial {\Bbb H}_{\bold C}^{n+1}=\{ (z, y_{n+1})
  \in {\bold C}^n \times {\bold R}: y_{n+1}-|z|^2=0 \}
  \times \{ t: t \in {\bold R} \}.$$
The first one is called the constraint boundary, the
second is called the free boundary. What we need is
the second one. Set ${\Cal H}:=\{ t+i \rho: t \in
{\bold R}, \rho>0 \}$ and let $\Omega(\Gamma,
{\Cal H})$ be the region of discontinuity of $\Gamma
\leq G({\bold Z})$ corresponding to ${\Cal H}$.
For $Z=(z, t, \rho)=(z, z_{n+1}) \in {\Bbb H}_{\bold
C}^{n+1}$ and $W=(w, \zeta, 0)=(w, w_{n+1}) \in \partial
{\Bbb H}_{\bold C}^{n+1}$, where $\rho(W)=\text{Im}
w_{n+1}-|w|^2=0$. If $w=0$, $\rho(W)=0$ implies that
$\text{Im} w_{n+1}=|w|^2=0$, i.e., $w_{n+1}=\zeta \in {\bold R}$.
The Poisson kernel associated with the free boundary is:
$$P(Z, \zeta):=\frac{\rho}{(\rho+\beta)^2+(t-\zeta)^2}.\eqno (6.4)$$

\proclaim\nofrills{\smc Theorem 6.1}. {\it $L P(Z, \zeta)^s=s(s-n-1)
P(Z, \zeta)^s.$}
\endproclaim

\proclaim\nofrills{\smc Theorem 6.2}. {\it $P(\gamma(Z), \gamma(\zeta))
|\gamma^{\prime}(\zeta)|=P(Z, \zeta)$, for $\gamma \in G({\bold Z})$.}
\endproclaim

  We define the Eisenstein series
$$E(Z, \zeta; s):=\sum_{\gamma \in \Gamma}
  P(\gamma(Z), \zeta)^{s}, \quad \text{Re}(s)>
  \delta(\Gamma), \eqno (6.5)$$
where $\delta(\Gamma)$ is the critical exponent of $\Gamma$.

\proclaim\nofrills{\smc Theorem 6.3}. {\it The Eisenstein series satisfies
the following properties:
$$E(\gamma(Z), \zeta; s)=E(Z, \zeta; s), \quad
  E(Z, \gamma(\zeta); s)=|\gamma^{\prime}(\zeta)|^{-s}
  E(Z, \zeta; s),$$ 
$$L E(Z, \zeta; s)=s(s-n-1) E(Z, \zeta; s). \eqno (6.6)$$}
\endproclaim

  The scattering matrix is defined as
$$S(\zeta, \eta; s):=\sum_{\gamma \in \Gamma}
  \frac{|\gamma^{\prime}(\zeta)|^s}
  {|\gamma(\zeta)-\eta|^{2s}}, \quad
  \text{Re}(s)>\delta(\Gamma),\eqno (6.7)$$
where $\zeta, \eta \in \Omega(\Gamma, {\Cal H})$.
It describes the normalized free boundary behaviour
of the Eisenstein series
$$S(\zeta, t; s)=\lim_{\rho \to 0} \lim_{\beta
  \to 0} \rho^{-s} E(Z, \zeta; s).$$

  For $\text{Re}(s)>n+1$, we define
$$G_0(Z, Z^{\prime}; s):=r_s(u(Z, Z^{\prime})), \quad
  \text{and} \quad G(Z, Z^{\prime}; s):=\sum_{\gamma \in
  \Gamma} G_0(\gamma(Z), Z^{\prime}; s),$$
where $r_s(u)=g_1(x)$ in Theorem 5.3.
Then, we have
$$\lim_{\rho^{\prime} \to 0} \lim_{\beta^{\prime} \to 0}
  (\rho^{\prime})^{-s} G_0(Z, Z^{\prime}; s)=c(s)
  P(Z, t^{\prime})^s.$$
Consequently,
$$\lim_{\rho^{\prime} \to 0} \lim_{\beta^{\prime} \to 0}
  (\rho^{\prime})^{-s} G(Z, Z^{\prime}; s)=c(s)
  E(Z, t^{\prime}; s).$$

\vskip 0.5 cm
\centerline{\bf 7. The modular forms and modular varieties
on $U(n+1, 1)$}
\vskip 0.5 cm

\flushpar
  In his paper \cite{10}, Wirthm\"{u}ller gave the Jacobi
modular forms associated to the root systems. Now we give
the definition of modular forms on $U(n+1, 1)$ associated
to $G({\bold Z})$.

\definition{\it Definition 7.1} {\it A modular form on $U(n+1, 1)$ associated
with $G({\bold Z})$} is a function
$\phi: {\Bbb H}_{\bold C}^{n+1} \to {\bold C}$ satisfying
the following transform equations:
\roster
\item
$\phi(\frac{z}{cz_{n+1}+d}, \frac{az_{n+1}+b}{cz_{n+1}+d})
=(c z_{n+1}+d)^k e^{2 \pi i m c(z_1^2+\cdots+z_n^2)/(c z_{n+1}+d)}
\phi(z, z_{n+1}).$
\item
$\phi(wz, z_{n+1})=\phi(z, z_{n+1})$ for all $w \in S_n$, where
$S_n$ is the symmetric group of $n$-order.
\item
$\phi(z, z_{n+1})$ is a locally bounded function as
$\text{Im} z_{n+1} \to \infty$.
\endroster
\enddefinition

  According to \cite{10}, for $Z=(z, z_{n+1}) \in
{\Bbb H}_{\bold C}^{n+1}$ and $w \in {\bold C}$, we have
$$\gamma(z, z_{n+1}, w)=(\frac{z}{cz_{n+1}+d}, \frac{az_{n+1}
  +b}{cz_{n+1}+d}, w-c \frac{z_1^2+\cdots+z_n^2}{cz_{n+1}+d}).$$

  The modular forms on $U(n+1, 1)$ can be written as follows
$$\phi \circ \gamma(z, z_{n+1})=(\frac{d \gamma(z_{n+1})}{dz_{n+1}})
  ^{-\frac{k}{2}} e^{-2 \pi im (\gamma(w)-w)} \phi(z, z_{n+1}).
  \eqno (7.1)$$

  Set ${\Cal M}_{n+1}=G({\bold Z}) \backslash
{\Bbb H}_{\bold C}^{n+1}$. It is well known that
${\Cal M}_1 \cong {\bold C}$, the coarse
moduli space of complex elliptic curves. 

  Now, we consider the modular functions on $U(n+1, 1)$.
$$(\phi|_{k, m} \gamma)(z, z_{n+1}) :=(cz_{n+1}+d)^{-k}
  e^m(\frac{-c(z_1^2+\cdots+z_n^2)}{cz_{n+1}+d}) \phi
  \circ \gamma(z, z_{n+1}), \eqno (7.2)$$
for $\gamma \in G({\bold Z})$, where $e^m(x):=e^{2 \pi
imx}$.

  The Eisenstein series for $U(n+1, 1)$ with weight $k$
and index $m$ is defined as
$$E_{k, m}(z, z_{n+1}) :=\sum_{\gamma \in G({\bold Z})_{
  \infty} \backslash G({\bold Z})}
  (1|_{k, m} \gamma)(z, z_{n+1}). \eqno (7.3)$$
Explicitly, this is
$$E_{k, m}(z, z_{n+1})=\frac{1}{2} \sum_{c, d \in {\bold Z},
  (c, d)=1} (cz_{n+1}+d)^{-k} e^m(\frac{-c(z_1^2+\cdots+z_n^2)}
  {cz_{n+1}+d}).$$
It is clear that
$$E_{k, m}(z, z_{n+1}+1)=E_{k, m}(z, z_{n+1}).$$
$$E_{k, m}(\frac{z}{z_{n+1}}, -\frac{1}{z_{n+1}})=\frac{1}{2}
  z_{n+1}^k \sum_{c, d \in {\bold Z}, (c, d)=1} (dz_{n+1}-c)^{-k}
  e^m(\frac{-c(z_1^2+\cdots+z_n^2)}{-cz_{n+1}+dz_{n+1}^2}).$$
By the identity:
$$\frac{-c}{-cz_{n+1}+dz_{n+1}^2}-\frac{1}{z_{n+1}}=\frac{-d}
  {dz_{n+1}-c},$$
we have
$$E_{k, m}(\frac{z}{z_{n+1}}, -\frac{1}{z_{n+1}})=z_{n+1}^k
  e^m(\frac{z_1^2+\cdots+z_n^2}{z_{n+1}}) E_{k, m}(z, z_{n+1}).$$

  The set of modular forms on $U(n+1, 1)$ is denoted as
$M_{k, m}(G({\bold Z}))$. It is obvious that
$E_{k, m}(z, z_{n+1}) \in M_{k, m}(G({\bold Z}))$.

  Now, let us define a family of modular functions for
$U(n+1, 1)$ of weight $0$ and index $0$, the $j$-invariants
associated to $U(n+1, 1)$:
$$j_m(z, z_{n+1}) :=\frac{1728 g_{2, m}(z, z_{n+1})^3}
 {\Delta_m(z, z_{n+1})}, \eqno (7.4)$$
where
$g_{2, m_1}(z, z_{n+1}) :=\frac{4}{3} \pi^4
 E_{4, m_1}(z, z_{n+1}),$
$g_{3, m_2}(z, z_{n+1}) :=\frac{8}{27} \pi^6
 E_{6, m_2}(z, z_{n+1}),$
and
$\Delta_m(z, z_{n+1}) :=g_{2, m}(z, z_{n+1})^3-27
  g_{3, \frac{3}{2}m}(z, z_{n+1})^2.$
In fact, for $\gamma \in G({\bold Z})$,
$$j_m(\gamma(z, z_{n+1}))=j_m(z, z_{n+1}). \eqno (7.5)$$
Therefore, $j_m$ is a modular function on the modular
variety ${\Cal M}_{n+1}$.

  If $z=0$ then ${\Cal M}_{n+1}$ degenerates to ${\Cal M}_1$
and $j_m(0, z_{n+1})=j(z_{n+1})$. It is well known that for every
$c \in {\bold C}$, $j(z_{n+1})=c$ has exactly one solution.
Thus, $j(z_{n+1})$ is an analytic isomorphism from ${\Cal M}_{n+1}$
to ${\bold C}$. Therefore, $j_m: {\Cal M}_{n+1} \to {\bold C}$
is a surjective morphism. Now, we complete the proof of Main
Theorem.

\enddocument